\documentclass[twoside, 12pt]{report}
\usepackage[english,russian]{babel}
\usepackage{epsfig}
\usepackage{psfig}
\usepackage{amssymb,amsmath,amsfonts,soul,amsthm,enumerate}
\usepackage[cp1251]{inputenc}
\usepackage[T2A]{fontenc}
\usepackage{mathrsfs}
\newtheorem*{Proposition}{Proposition}

\def\@evenfoot{}
\def\@oddfoot{}

\begin{document}
\def\@evenhead{\vbox{\hbox to \textwidth{\thepage\leftmark}\strut\newline\hrule}}

\def\@oddhead{\raisebox{0pt}[\headheight][0pt]{%
\vbox{\hbox to \textwidth{\rightmark\thepage\strut}\hrule}}}

\newpage
\normalsize
\def\bibname{\vspace*{-30mm}{\centerline{\normalsize References}}}
\thispagestyle{empty}
\begin{flushleft}
{\small {\bf EURASIAN  JOURNAL OF MATHEMATICAL \\ AND COMPUTER APPLICATIONS}\\
 ISSN 2306--6172\\
Volume 4, Number   (3)    --}
\end{flushleft}
\vskip 5 mm

\centerline{\bf INVERSE SOURCE PROBLEM FOR WAVE EQUATION }
\centerline{\bf AND GPR DATA INTERPRETATION PROBLEM}

\vskip 0.3cm
\centerline{\bf Balgaisha Mukanova, Vladimir G. Romanov}

\vskip 0.3cm
\noindent{\bf Abstract}
{\small The inverse problem of identifying the unknown spacewise dependent source $F(x)$ in $1D$ wave equation  $u_{tt}=c^2u_{xx}+F(x)H(t-x/c)$, $(x,t)\in \{(x,t)| x>0, -\infty\le t\le T\}$ is considered. Measured data are taken in the form $g(t):=u(0,t)$.
The relationship between that problem and Ground Penetrating Radar (GRR) data interpretation problem is shown. The non-iterative algorithm for reconstructing the unknown  source $F(x)$ is developed. The algorithm is based on the Fourier expansion of the source $F(x)$ and the explicit representation of the direct problem solution via the function $F(x)$. Then the minimization problem for discrete form of the Tikhonov functional is reduced to the linear algebraic system and solved numerically.
Calculations show that the proposed algorithm allows to reconstruct the spacewise dependent source $F(x)$ with enough accuracy for noise free and noisy data.
}
\vskip 0.3cm

\setcounter{figure}{0}

\renewcommand{\thesection}{\large 1}

\section{\large Introduction}

In this paper we study the problem of identifying an unknown spacewise dependent source $F(x)$ in
\begin{eqnarray}\label{dp1}
\left\{\begin{array}{ll}
u_{tt} -c^2u_{xx}=F(x) H(t-x/c),~c=const>0,
\\(x, t)\in\Omega_T=\{(x,t)~|~ x>0, \> -\infty\le t\le T\};
\\ \left( u_t -c_0u_x\right)_{x=0}=0,~c_0=const>0, \quad
\quad u|_{t<0}=0,
\end{array}
\right.
\end{eqnarray}
from boundary measured data
\begin{eqnarray}\label{measured data}
g(t): = u(0,t), \quad t\in [0,T].
\end{eqnarray}
Here the function $F(x)$ is assumed to have a finite support in $(0, \infty)$ and $H(t)$ is a given piecewise smooth function such that $H(t)\equiv 0$ for $t<0$ and $H(+0)\neq0$. We define this problem as an inverse source problem (ISP) for wave equation \eqref{dp1} with Dirichlet type boundary measured data (\ref{measured data}).

Inverse problems for hyperbolic equations naturally arise from medical applications, seismology and geophysical prospecting, radar technology, electrical networks  and many other physical problems (see \cite{cannon}-\cite{Tikh} and references therein). An inverse source problem of identifying an unknown source term $S(u)$ in the wave equation $u_{tt}-u_{xx}=S(u)$, $x,t > 0$, from boundary data $u(0,t)=f(t)$, $u_x(0,t)=g(t)$, has first been studied in \cite{cannon}. Here an existence result for the identification problem is derived.
Unicity of the solution and ill-conditionedness of the ISP for wave equation  with variable speed function and final measured data has been studied in \cite{AH1}.
Uniqueness results for multidimensional parabolic and hyperbolic ISPs have been established in \cite{Engl1}. Stability estimate and a reconstruction formula for $f(x)$ in the hyperbolic equation $u_{tt}=\Delta u+\sigma (t)f(x)$, $x\in \Omega \subset \mathbb{R}^r$, $t>0$, from the Neumann type additional data $\partial u(x,t;f)/\partial n$ have been obtained in \cite{Yamam}.
Regarding the numerical approaches to hyperbolic coefficient inverse problems, we refer to monographs \cite{Kab2}, \cite{Klib1}.

Most of numerical approaches to ISP for parabolic and hyperbolic equations deal with source term in separable form $F(x)H(t)$ (see, for instance, \cite{chapouly}, \cite{Yamam}, \cite{AH2}, \cite{Sebu} and references therein). In this paper the function $H(x,t)$ has the form $H(t-x/c)$, since, as it is shown below, the linearized GPR data interpretation problem takes the form (\ref{dp1})-(\ref{measured data}); therefore the proposed method is applicable in radar techniques. In practice inverse problems arising in GPR techniques are solved via different approximate ways, most relevant of them are described in \cite{Stolt} and in the review \cite{Caner}.

In this paper, we develop new non-iterative algorithm for identifying the spacewise dependent source $F(x)$ in (\ref{dp1})-(\ref{measured data}).  This algorithm is based on integral formula for the solution of wave equation (\ref{dp1}) and use of the $N$th partial sum of the Fourier expansion for the term $F(x)$. Substituting then this formula in the regularized cost functional
\begin{eqnarray}\label{functional-1}
J_{\alpha}(F): = \frac{1}{2}\Vert u(0,\cdot;F)-g(\cdot) \Vert_{L^2(0,T)}^2+\frac{\alpha}{2} \Vert F \Vert_{L^2(0,l)}^2, \quad \alpha>0,
\end{eqnarray}
where $l=l(T)$,  we obtain a system of algebraic equations which unique solution gives an approximate regularized solution of the considered inverse problem. The algorithm is simple, effective and does not require any iterative procedures. Our numerical results demonstrate that the accuracy of all reconstructions are sufficient for high noise levels of measured data. The similar approach for an inverse source problem related to the advection--diffusion equation has been proposed in \cite{AH2}, \cite{Sebu}.

The paper is organized as follows. In Section 2 we reduce the GPR data interpretation problem to the ISP (\ref{dp1})-(\ref{measured data}). Numerical algorithm for identification of a spacewise dependent source
 from Dirichlet type measured output data is described in Section 3. Results of computational experiments
 are given in Section 4. Some concluding remarks are made in Section 5.

\renewcommand{\thesection}{\large 2}
\section{\large Linearized  mathematical model of GPR method }
  Let us formulate the 1D inverse problem for the model of GPR technique. As it is common in geophysics, assume that the medium fills the half-space $z>0$ and the half-space $z<0$ corresponds to the air. Let the electrical permittivity $\varepsilon$ of the medium depend on the coordinate $z$ only, magnetic permittivity  $\mu=\mu_0=\mbox{const}>0$ in the whole space and the conductivity is negligible. Let the current source with intensity
\begin{eqnarray*}
  j^{ex}(t)=\Phi(t)\delta(z),~ \Phi(t)=0 \mbox{ if } t\leq 0,~ \Phi(t)\in C^2[0,  \infty),~\Phi''(+0)\neq 0,
\end{eqnarray*}
be placed at the boundary $z=0$ and directed along the axis $y$.
Then it follows from Maxwell's equations that the electromagnetic field depends on $(z,t)$ only. The field has an electric component $E_2(z,t)$ along the axis $y$, and a magnetic component $H_1(z,t)$  along the axis $x$ that satisfy the Cauchy problem:
\begin{eqnarray}\label{Maxwell}
\frac{\partial H_1}{\partial z} =\varepsilon(z) \frac{\partial E_2}{\partial t} +\delta(z) \Phi(t),
\quad
\frac{\partial E_2}{\partial z} =\mu(z) \frac{\partial H_1}{\partial t},
\quad  (E_2, H_1)_{t<0}=0.
\end{eqnarray}

 We assume below $\mu(z)=\mu_0>0$.
By taking first derivatives with respect to $t$ from first equation and with respect to $z$ from second one in (\ref{Maxwell}) and eliminating $\partial^2 H_1/\partial t \partial z$ we get

\begin{eqnarray}\label{dpGPR}
\frac{\partial^2 E_2}{\partial z^2}
=\mu_0\varepsilon(z) \frac{\partial^2 E_2}{\partial t^2} +\mu_0\delta(z) \Phi'(t),
\quad
\quad E_2|_{t<0}=0.
\end{eqnarray}

Denote by $c(z)=1/{\sqrt{\mu_0\varepsilon(z)}}$. Suppose that the function $c^{-2}(z)$ is presented in the following form
\begin{eqnarray}\label{cz}
c^{-2}(z)=\left\{\begin{array} {ll}c_0^{-2},& \mbox{if } z<0\\
c_1^{-2}+F(z),& \mbox{if } z\geq 0,\end{array}\right.
\\
c_0, c_1=const,~ F(z)\in C(R),~ |F(z)|\ll  c_1^{-2},
\end{eqnarray}
where the function $F(z)$ has a finite support in $z\in (0,\infty)$ and values $c_0>0,~ c_1>0$ are given. As it has been shown in \cite{KabIsk}, the  conditions imposed to the functions $\Phi(t),~ c^2(z)$ provide the existence and uniqueness of the solution to the Cauchy problem (\ref{dpGPR}).

Now represent the solution of the direct problem (\ref{dpGPR}) in the form
$E_2(z,t)=U(z,t)+u(z,t)$ where $U(z,t)$ is the generalized solution of the Cauchy problem:
\begin{eqnarray}\label{dpU}
\begin{array}{ll}
U_{zz}=\frac{1}{\overline{c}^{2}(z)}U_{tt}+\mu_0\Phi'(t)\delta(z),\quad  (z\in \mathbb{R},~ t>-\infty), \\
U|_{t<0}\equiv0,\\
\overline{c}^2(z)=\left\{\begin{array} {ll}c_0^2,& \mbox{if } z<0,\\
c_1^2,& \mbox{if } z\geq 0.\end{array}\right.
\end{array}
\end{eqnarray}

Then the solution of the problem (\ref{dpU}) is given by the formula:
 \begin{eqnarray}\label{U}
U(z,t)=-\frac{\mu_0c_0c_1}{c_0+c_1}\left\{
\begin{array}{ll}
\Phi(t+z/c_0), & z<0,
\\
\Phi(t-z/c_1), & z>0.
\end{array}
\right.
\end{eqnarray}
It can be checked directly that the function $U(z,t)$ is continuous anywhere and twice continuously differentiable in the half spaces $\mathbb{R}^2_-=\{(z,t)| ~z<0, t\in\mathbb{R}\}$,  $\mathbb{R}^2_+=\{(z,t)| ~z>0, t\in\mathbb{R}\}$ and its first derivatives at $z=0$ are expressed as
 \begin{eqnarray*}
U_z(-0,t)=-\frac{\mu_0 c_1}{c_0+c_1}\Phi'(t), \quad U_z(+0,t)=\frac{\mu_0 c_0}{c_0+c_1}\Phi'(t),
\end{eqnarray*}
i.e.
 \begin{eqnarray*}
U_z(+0,t)-U_z(-0,t)=\mu_0\Phi'(t).
\end{eqnarray*}
 The last formula confirms that the second derivative $U_{zz}$ is represented as the singular function $\mu_0\Phi'(t)\delta(z)$ and a regular one.

 The linearization of the equation (\ref{dpGPR}) with respect to $u(z,t)$ shows that the function $u(z,t)$ satisfies the equation
 \begin{eqnarray*}
\frac{\partial^2 u}{\partial z^2}
=\frac{1}{\bar{c}^2(z)} \frac{\partial^2 u}{\partial t^2} +F(z) \frac{\partial^2 U}{\partial t^2},
\quad
\quad u|_{t<0}=0.
\end{eqnarray*}

Since the support of the function $F(z)$ belongs to the domain  $z>0$, the function  $u(z,t)$ with its first derivatives are continuous at the axis $z=0$. For  $z<0$  the function $u(z,t)$ is a solution of the homogeneous  equation
  and is expressed in the form  $u(z,t)=r(t+z/c_0)$, where $r(t)=u(0,t)$. Therefore it satisfies the condition
$ u_t-c_0u_z=0$ for $z\le 0$ and, by the continuity, for $z=+0$. Then for $z>0$ the function $u(z,t)$ is a solution to the following problem
 \begin{eqnarray}\label{dp3}
c_1^2\frac{\partial^2 u}{\partial z^2}
=\frac{\partial^2 u}{\partial t^2} +F(z) c_1^2\frac{\partial^2 U}{\partial t^2}, \>z>0;
\quad \left(\frac{\partial u}{\partial t} -c_0\frac{\partial u}{\partial z}\right)_{z=0}=0, \quad
\quad u|_{t<0}=0.
\end{eqnarray}

In GPR method electrical field $E_2(0,t)$ is measured, therefore  additional data for inverse problem are
\begin{equation}\label{adddata}
u|_{z=0}=g(t)\equiv E_2|_{z=0}-U|_{z=0}, \quad t\in[0, T], ~ T>0.
\end{equation}
Now introduce the notation
\begin{equation}\label{Ht}
H(t)=\frac{\mu_0c_0c_1^3}{c_0+c_1}\Phi''(t),
\end{equation}
replace $z$ by $x$ and $c_1$ by $c$ in (\ref{dp3})  and define $\Omega_T=\{(x,t)| x>0, -\infty\le t\le T\}$.   Then the direct problem for $u(x,t)$   is formulated as follows
 \begin{eqnarray*}
\frac{\partial^2 u}{\partial t^2} -c^2\frac{\partial^2 u}{\partial x^2}=F(x) H(t-x/c),\quad (x, t)\in\Omega_T;
\\
 \left(\frac{\partial u}{\partial t} -c_0\frac{\partial u}{\partial x}\right)_{x=0}=0, \quad
\quad u|_{t<0}=0,
\end{eqnarray*}
which coincides with the direct problem statement (\ref{dp1}).

Therefore, the GPR data interpretation problem  is reduced to the linear ISP (\ref{dp1})-(\ref{measured data}).

Note that $H(t)=0$ for $t<0$, then $u(x,t)\equiv 0$ for $t\le x/c$.
Moreover, to calculate $g(t)=u(0,t)$ for $t\in[0,T]$
we need only to find the solution of (\ref{dp1}) in the domain
\begin{equation}\label{DomainD}
D_T=\{(x,t)~| ~0\le x/c\le t\le T-x/c\}.
\end{equation}

\medskip
{\Proposition If $H(t)\in H^1[0,T]$, $ H(0)\neq0$ and  $ g(t)\in H^2[0,T]$ then the space-dependent source $F(x)\in L^2(0,l)$ for $x\in [0, l]$, $l= cT/2$, for ISP (\ref{dp1}), can be identified uniquely from the boundary measured data (\ref{measured data}).}

\medskip
\noindent{\textbf{Proof}}
By introducing the notation
\begin{equation}
v(x,t)=\left(\frac{\partial u}{\partial t} +c\frac{\partial u}{\partial x}\right)
\end{equation}
the equation (\ref{dp1}) is rewritten as follows:

\begin{equation}\label{eq_v}
\frac{\partial v}{\partial t} -c\frac{\partial v}{\partial x}=F(x)H(t-x/c).
\end{equation}
Since $u(x,t)=0$, then $v(x,t)=0$ for $\{(x,t)~|~ 0\le t\le x/c\}$, and, in particulary, $v(x, x/c)=0$. Let $(x_0,t_0)$ be an arbitrary point in $D_T$. Integrate the equation (\ref{eq_v}) along  the line $t+x/c=t_0+ x_0/c$ from the point $(x_0,t_0)$ up to the intersection with the characteristic line $t=x/c$, i.e. the point $((x_0+ct_0)/2, (t_0+ x_0/c)/2)$, and obtain

\begin{equation*}
 v(x_0,t_0)=\frac{1}{c}\int\limits_{x_0}^{(x_0+ct_0)/2} F(x) H(t_0+ x_0/c-2x/c) dx, \quad (x_0,t_0)\in D_T.
\end{equation*}
Then for all $(x,t)\in D_T$ the following equality
\begin{equation}\label{utux}
\left( \frac{\partial}{\partial t} +c\frac{\partial }{\partial x}\right)u(x,t)=\frac{1}{c}\int\limits_{x}^{(x+ct)/2} F(\xi) H(t+ x/c-2\xi/c) d\xi, \> (x,t)\in D_T.
\end{equation}
holds. The combination of the expression above at $x=0$ with the boundary condition
\begin{equation*}
\left(\frac{\partial u}{\partial t} -c_0\frac{\partial u}{\partial x}\right)_{x=0}=0
\end{equation*}
defines the derivative
\begin{equation}\label{ut}
\left. \frac{\partial u}{\partial t}\right|_{x=0}=\frac{c_0}{c(c+c_0)}\int\limits_{0}^{ct/2} F(\xi) H(t-2\xi/c) d\xi.
\end{equation}

 It follows from expressions (\ref{ut}) that
 \begin{equation}\label{g1}
g'(t)=\frac{c_0}{c(c+c_0)}\int\limits_{0}^{ct/2} F(\xi) H(t-2\xi/c) d\xi, \quad t\in[0,T]
 \end{equation}
 Let us analyze the expression  (\ref{g1}). If $F(x)\in L^2(0,l)$, $l=T/(2c)$, $H(t)\in H^1[0,T]$ then  $g(t)\in H^2[0,T]$ and $g'(0)=0$.  Taking first derivative of  (\ref{g1}), we have
\begin{equation}\label{g2}
g''(t)=\frac{c_0}{c(c+c_0)}\Big[F(ct/2)H(0)c/2+
\int\limits_{0}^{ct/2} F(\xi) H'(t-\xi/(2c)) d\xi\Big], \quad t\in[0,T].
\end{equation}
If $H(0)\neq0$ then (\ref{g1}) is rewritten as follows:
\begin{equation}\label{g3}
\hat{g}(x)=F(x)+
\frac{2}{c H(0)}\int\limits_{0}^{x} F(\xi) H'(2(x-\xi)/c) d\xi, \quad x\in[0,l],
\end{equation}
where $\hat{g}(x)=2(c+c_0)g''(2x/c)/(c_0 H(0))$ and $l=cT/2$.
The equation (\ref{g3}) represents Volterra equation of the second kind and is uniquely solvable in $ L^2(0,l)$  for all $\hat{g}(x)\in L^2(0,l)$ (\cite{Tricomi}). In other words, boundary data $g(t)$, $t\in[0,T]$  uniquely define the function $F(x)$ for $x\in[0,l], l=cT/2$.
\hfill$\Box$

\textbf{Remark} The equation (\ref{g3}) gives an alternate way of solving the ISP (1)-(2). For instance, it can be solved numerically. The bigger values of $|H(0)|$ correspond to better stability estimates for the solution of the equation (\ref{g3}). This observation is in the concordance with numerical results presented below.

\renewcommand{\thesection}{\large 3}
\section{\large Algorithm for identifying the spacewise dependent source}

Now we are going to construct a computational algorithm for solving the considered ISP.  For a given $F \in L^2(0,l)$  denote by $u:=u(x,t;F)$ a solution of the direct problem (\ref{dp1}).
Derive the formula for that solution at the axis $x=0$.
It follows from (\ref{ut}) and initial condition $u(0,0)=0$ that
 \begin{equation}\label{u}
  u(0,t;F)=\frac{c_0}{c(c+c_0)}\int\limits_{0}^{t}\int\limits_{0}^{c\tau/2} F(\xi) H(\tau-2\xi/c) d\xi d\tau.
 \end{equation}

We assume now that the function $F(x)$ has a finite support in $(0,l), l=cT/2$,  and approximate the unknown source $F(x)$ by the $N$th partial sum of the Fourier series at the interval $[0,l]$:
\begin{eqnarray}\label{FFourier}
F^N(x)=\sum\limits_{k=1}^{N}F_k X_k(x),
\end{eqnarray}
where $X_k(x), ~k=\overline{1,\infty}$, are eigenfunctions of the following spectral problems:
\begin{eqnarray}\label{sp}
\left \{ \begin{array}{ll}
X_k''+\lambda_k^2X_k=0,  \quad x \in (0,l);\\
X_k(0)=0, ~~X_k(l) = 0.
\end{array} \right.
\end{eqnarray}

Solving the two-point problem  (\ref{sp}), we find the normalized eigenfunctions
\begin{eqnarray}\label{X_k}
X_k(x)=\sqrt{\frac{2}{l}} \sin(\lambda_k x),~ \lambda_k=\frac{k\pi}{l},~ k=\overline{1,\infty},
\end{eqnarray}
corresponding to the eigenvalues $\lambda_k$.
Note that the eigenfunctions system $X_k(x)$,  $k=\overline{1,\infty},$ is complete in $L^2(0,l)$.

Substituting (\ref{FFourier}) into (\ref{u}) we get
\begin{eqnarray}\label{uDFN}
u(0,t;F^N)=\frac{c_0}{c(c+c_0)}\sum\limits_{k=1}^{N}F_k
\int\limits_{0}^{t}\int\limits_{0}^{c\tau/2} X_k(\xi) H(\tau-2\xi/c) d\xi d\tau.
\end{eqnarray}
Changing the integration order and introducing new variable $s=\tau-2\xi/c$ in (\ref{uDFN}) we obtain
\begin{eqnarray}\label{uFN}
u(0,t;F^N)=\frac{c_0}{c(c+c_0)}\sum\limits_{k=1}^{N}F_k \int_0^{ct/2}X_k(\xi)
\int_{0}^{t-2\xi/c}H(s)ds d\xi
\end{eqnarray}
By notation (\ref{Ht})
\begin{eqnarray}\label{Psi}
\int\limits_0^tH(s)ds=\frac{\mu_0 c_0 c^3}{c_0+c}(\Phi'(t)-\Phi'(0)).
\end{eqnarray}
Taking into account (\ref{uFN}) and  (\ref{Psi}) we have
\begin{eqnarray*}
\begin{array}{r}
u(0,t;F^N)=\frac{\mu_0}{(c^{-1}+c_0^{-1})^2}\sum\limits_{k=1}^{N}F_k \int\limits_0^{ct/2} X_k(\xi)(\Phi'(t-2\xi/c)-\Phi'(0)) d\xi = \sum\limits_{k=1}^{N}F_k G_k(t),
\end{array}
\end{eqnarray*}
where the following notation is used:
\begin{eqnarray}\label{uFN1}
\begin{array}{r}
G_k(t)\triangleq \frac{\mu_0}{(c^{-1}+c_0^{-1})^2}\int\limits_0^{ct/2} X_k(\xi)(\Phi'(t-2\xi/c)-\Phi'(0))d\xi,
~k=\overline{1,N}. \end{array}
\end{eqnarray}

Because the measured data $g(t)$ always contain a random noise, we look for the unique regularized solution of the inverse problem (\ref{dp1})-(\ref{measured data}). This solution $F_{\alpha} \in L^2(0,l)$ is defined as a minimum of the Tikhonov functional (\ref{functional-1}).
The regularized cost functional (\ref{functional-1}) on the finite-dimensional approximation $F^N(x)$ is the  $N$-variable function $J_{\alpha}(F^N) \equiv J_{\alpha}(F_1^N,F_2^N, \cdots F_N^N)$:
\begin{eqnarray*}
J_{\alpha}(F^N)\triangleq \frac{1}{2}\int\limits_{0}^{T} \left(\sum \limits_{k=1}^{N} F^N_k G_k(t)-g(t)\right)^2dt+\frac{\alpha}{2}\sum \limits_{k=1}^{N}\left(F_k^N \right)^2.
\end{eqnarray*}
The $N$-dimensional vector of unknown parameters $(F_1^N,F_2^N,\ldots, F_N^N)$ is the unique minimizer of this functional and is defined from the conditions
\begin{eqnarray*}
\begin{array}{r}
\frac{\partial J_{\alpha}(F_1^N,F_2^N,\ldots, F_N^N)}{\partial F_k^N}:=
\sum \limits_{i=1}^{N} F_i^N \int\limits_{0}^{T}G_i(t)G_k(t) dt +\alpha F_k^N
- \int \limits_{0}^{T} G_k(t)g(t)dt=0, ~ k=\overline{1,N}.
\end{array}
\end{eqnarray*}
This yields the following system of linear algebraic equations
\begin{eqnarray}\label{linear_system}
({\mathbf A^N}+\alpha{\mathbf I}){\mathbf F^N_{\alpha}} ={\mathbf b^N},
\end{eqnarray}
with respect to the unknown vector $\mathbf F^N_{\alpha}:=(F_{\alpha 1}^N,F_{\alpha2}^N,\ldots, F_{\alpha N}^N)$, with the matrix ${\mathbf A^N}$ and right
hand side vector ${\mathbf b^N}$, defined as
\begin{eqnarray}\label{c16}
\begin{array}{ll}
{ A^N}_{ij}=\int\limits_{0}^{T}G_i(t)G_j(t)dt,~~i,j=\overline{1,N},\\
{ b^N}_j=\int \limits_{0}^{T}G_j(t)g(t)dt ,~~j=\overline{1,N}.
\end{array}
\end{eqnarray}
Here ${\mathbf I}$ is the identity matrix and the functions $G_j(t)$  are defined by (\ref{uFN1}).
Hence, the unique solution of the discrete problem (\ref{linear_system})-(\ref{c16}) defines an approximate solution of the regularized inverse problem.
The problem of choosing the regularization parameters $N$ and $\alpha$ will be discussed in the next section.

\renewcommand{\thesection}{\large 4}
\section{\large Numerical results}
{
Before using the described algorithm, let us analyze the behavior of the relative error for different values of the parameter of regularization $\alpha>0$, cut-off parameter $N$ and noise level $\gamma>0$.

 Let
\begin{eqnarray*}
\varepsilon_F :=\Vert F-F^N \Vert_{L^2(0,1)}/\Vert F \Vert_{L^2(0,1)}.
\end{eqnarray*}

In order to obtain noise free synthetic measured data we have used the formula (\ref{u}) and calculated the integral in (\ref{u}) numerically. But in practice, measured data always contain noise, so, we define the random noisy output data as follows:
\begin{eqnarray*}
g^{\gamma}(t)=g(t)+\delta g(t)=g(t)+\gamma n(t)\|g(t)\|_{L^2[0,T]}/\|n(t)\|_{L^2[0,T]},
\end{eqnarray*}
where $\gamma>0$ is the relative noise level and
 \begin{eqnarray*}
 n(t)=\sum\limits_{j=0}^{N_n}\xi_j\eta\big(\frac{t-j\tau}{\tau}\big), \ \tau=T/N_n
\end{eqnarray*}
is the random function. Here $\eta(t)$ is a standard linear finite element and the values  $\xi_j$, $j=0,..,N_n$ are obtained using the MATLAB "randn" function, which generates arrays of random numbers whose elements are normally distributed with mean $0$ and standard deviation $\sigma=1$.

Let us assume now that the right hand side of the linear system (\ref{linear_system}) contains an error $\delta \mathbf{b}^N$. Then the relative error of the solution $\delta \mathbf{F}^N$, which is defined as the difference between solutions obtained for noise free and noisy data, is estimated as follows:
\begin{eqnarray}\label{cond}
|\delta \mathbf{F}^N|\leq  C(\mathbf{A}^N,\alpha) |\delta\mathbf{b}^N|,
\end{eqnarray}
here  $C(\mathbf{A}^N,\alpha)$  is a condition number of the matrix  $\mathbf{A}^N+\alpha \mathbf{I}$,  which depends on $N$, $\alpha$, $T$, $c$, $c_0$ and the function $H(\cdot)$ as well.
The expressions (\ref{c16}) show that the errors in coordinates of $\delta \mathbf{b}_N$ in (\ref{cond}) are estimated via the relative noise level $\gamma=\|\delta g(t)\|_{L_2}/\|g(t)\|_{L_2}$ as follows:
\begin{eqnarray}\label{rhs_error}
|\delta b_j^N|=\left |\int \limits_{0}^{T} G_j(t)\delta g(t)dt\right |\leq C_1\|\delta g(t)\|_{L_2(0,T)}
= C_1\gamma\| g(t)\|_{L_2(0,T)}, 
\end{eqnarray}
where $C_1=\max\limits_{1\leq j\leq N}\|G_j(t)\|_{L_2(0,T)}$.
This yields:
\begin{eqnarray}\label{rel_err_noise}
|\delta \mathbf{F}^N|\leq C(\mathbf{A}^N,\alpha) \sqrt{\sum\limits_{j=1}^N (\delta \mathbf{b}_j^N)^2}
\leq C(\mathbf{A}^N,\alpha)\sqrt{N}C_1 \gamma\| g(t)\|_{L_2(0,T)}.
\end{eqnarray}
Therefore, the estimate (\ref{rel_err_noise}) establishes relationship between relative error of the approximate solution $F^N(x)$ for noisy data and the noise level $\gamma$. This estimate also shows that the most admissible parameters $N$ and $\alpha$ should correspond to minimal value of the number $\sqrt{N}C(\mathbf{A}^N,\alpha)$.  The last point leads to the practical way to choose these parameters.

The additional analysis has been done by computing the values of discrepancy $\eta$ defined as follows
\begin{equation}\label{eta}
\eta=\left(\int\limits_{0}^{T} \left(\sum \limits_{k=1}^{N} F^N_{\alpha k} G_k(t)-g(t)\right)^2dt\right)^{1/2}.
\end{equation}

Let the assumptions of the Proposition hold and $F(x)$ be the exact solution of the considered ISP. Let $F^N(x)$ and $F^N_{ex}(x)$ be computed and exact versions of the partial Fourier sums of $F(x)$. Denote by $C_i, i=1,2,3$ different constants  which do not depend on $F(x)$ and can depend on $N$, $\alpha$ and physical parameters of the problem.  Then the difference between exact and numerical solution of the inverse problem is estimated as follows:
\begin{equation}\label{deltaF}
  \|F(x)-F^N(x)\|_{L^2(0,l)} \leq \|F(x)-{F}_{ex}^N(x)\|_{L^2(0,l)}+ \|F_{ex}^N(x)-F^N(x)\|_{L^2(0,l)}.
\end{equation}

Define the function $g^N(t)=u(0,t;F_{ex}^N)$. Subtracting the equation (\ref{uFN}) from (\ref{u}) we obtain the integral equation which links the functions
 $g(t)- g^N(t)$  and $F(x)-F_{ex}^N(x)$. The solution of that equation satisfies the stability estimate which can be obtained in standard way:
\begin{equation}\label{deltaFN}
   \|F(x)-F_{ex}^N(x)\|_{L^2(0,l)} \leq C_2\|g(t)- g^N(t)\|_{H^1[0,T]}.
\end{equation}

Due to the orthogonality of basic functions $X_k(x)$  the $L_2$ -norm of the function
$\delta F^N(x)=F_{ex}^N(x)-F^N(x)$
 is equal to Euclidian norm of the vector $\delta \mathbf{F}^N(x)$; therefore combination of (\ref{deltaFN}) with  (\ref{rel_err_noise}) estimates the computational error of the solution to the inverse problem:
\begin{equation}\label{deltaFlast}
\|F(x)-F^N(x)\|_{L^2(0,l)} \leq C_2\|g(t)- g^N(t)\|_{H^1[0,T]}  + C_3 \|\delta g^N(t)\|_{L^2[0,T]}.
\end{equation}
}

As it is seen from the definition of the matrix $\mathbf{A}$, it can be calculated independently before measurements. Therefore the condition numbers  $C(\mathbf{A}^N,\alpha)$   for different values of $N$, $\alpha$ and given physical data  $c$, $c_0$, $T$, $H(t)$, $l$  can be defined. Then the most admissible combinations of $N$ and $\alpha$ can be established.
\begin{table}
{\small
\noindent\textbf{Table 1.} Values of the condition number  $C(\mathbf{A}^N,\alpha)$  depending on the parameters $N$, $\alpha$ for different $\Phi(t)$, $\beta_1=1.546$, $\beta_2=1.373$, $T=12\cdot10^{-9}$ sec, $c=1.5 \cdot 10^8$ m/sec, $l=0.9$ m:
\begin{center}
{\small \begin{tabular}{|c|c|c|c|c|c|c|c|c|c|c|c|}
\hline
\multicolumn{12}{|c|}{$~\Phi(t)=\sin(8t+\beta_1)\exp(-0.2t)~|~~\Phi(t)=\sin(t+\beta_2)\exp(-0.2t)$}
\\
 \hline $N \diagdown \alpha$ & $0$ &$10^{-4}$ &$10^{-3}$& $10^{-2}$ &$10^{-1}$ &$0$ &$10^{-5}$ &$10^{-4}$ &$10^{-3}$
 &$10^{-2}$ &$10^{-1}$
 \\ \hline
$5$ &$1.06$&$1.06$&$1.06$&$1.06$ & $1.057$ & $4.5$&$4.5$&$4.5$&$4.5$ &$4.45$ & $4.20$ \\ \hline
$8$ &$1.17$&$1.17$&$1.17$& $1.17$ & $1.16$ & $45.0$&$45.0$&$45$&$44.6$ &$41.4$ & $24.2$ \\ \hline
$11$ &$1.39$&$1.39$&$1.39$& $ 1.38$ & $1.36$ & $197$&$196.6$&$196$&$189$ &$142$ & $41.1$ \\ \hline
$14$ &$1.75$&$1.75$&$1.75$& $1.75$ & $1.71$ & $562$&$562$&$556$&$507$ &$268$ & $47.6$\\ \hline
$17$ &$2.43$&$2.43$&$2.43$& $ 2.42$ & $ 2.34$ & $1278$&$1275$&$1247$&$1022$ &$365$ & $50.1$ \\ \hline
$20$ &$3.72$&$3.72$&$3.72$& $3.70$ & $3.55$ & $2511$&$2498$&$2393$&$1685$ &$426$ & $51.1$ \\ \hline
\end{tabular}}
\end{center}}
\end{table}
Table 1 shows values of condition numbers  $C(\mathbf{A}^N,\alpha)$  computed for the function $\Phi(t)=\sin(\omega t+\beta)\exp(-\gamma t)-\Phi_0$  with different values of $\omega$, $\alpha$ and $N$. The parameters $\beta$ and $\Phi_0$ are taken to satisfy the conditions $\Phi(0)=0$, $\Phi'(0)=0$, namely, $\beta=\arctan(\omega/\gamma)$, $\Phi_0=\sin\beta$.
Values of discrepancies $\eta(N,\alpha)$ calculated for noise free data are collected in Table 2.

\begin{table}
{\small
\noindent\textbf{Table 2.}  Values of the discrepancy $\eta$ for noise free data depending on the parameters $N$, $\alpha$ and other inputs defined in Table 1:
\begin{center}
{\small \begin{tabular}{|c|c|c|c|c|c|c|}
\hline
\multicolumn{7}{|c|}{$~\Phi(t)=\sin(8t+\beta_1)\exp(-0.2t)$}
\\
 \hline $N \diagdown \alpha$ & $0$ &$10^{-5}$ &$10^{-4}$ &$10^{-3}$ &$10^{-2}$ &$10^{-1}$
 \\ \hline
$5$ &$0.084$&$0.084$&$0.084$&$0.084$&$0.084$
&$0.093$\\ \hline
$8$ &$0.054$&$0.054$&$0.054$&$0.054$&$0.054$&$0.067$\\ \hline
$11$ &$0.012$&$0.012$&$0.012$&$0.012$&$0.013$&$0.042$\\ \hline

$14$ &$0.0065$&$0.0065$&$0.0065$&$0.0065$ &$0.008$&$0.041$\\ \hline

$17$ &$0.0043$&$0.0043$&$0.0043$&$0.0043$
&$0.006$&$0.04$\\ \hline

$20$ &$0.0029$&$0.0029$&$0.0029$&$0.003$
&$0.005$&$0.04$\\ \hline

\multicolumn{7}{|c|}{$\Phi(t)=\sin(t+\beta_2)\exp(-0.2t)$}
\\
\hline $N \diagdown \alpha$ & $0$ &$10^{-5}$ &$10^{-4}$ &$10^{-3}$ &$10^{-2}$ &$10^{-1}$
\\ \hline
 $5$ &$0.0034$&$0.034$&$0.034$&$0.034$ &$0.034$&$0.049$\\ \hline
$8$ &$0.01$&$0.01$&$0.01$&$0.01$ &$0.011$&$0.039$\\ \hline
$11$ &$0.001$&$0.001$&$0.001$&$0.0011$ &$0.0043$ & $0.038$\\ \hline
$14$ &$0.00044$&$0.00044$&$0.00044$&$0.0006$ &$0.0043$ & $0.038$\\ \hline
$17$ &$0.00037$&$0.00036$&$0.00046$&$0.00054$ &$0.0043$ & $0.038$\\ \hline
$20$ &$0.00037$&$0.00036$
&$0.00036$&$0.00054$ &$0.0043$ & $0.038$\\ \hline
\end{tabular}}
\end{center}}
\end{table}

It is seen in Table 1 that the most important parameters that influence to the condition number are the frequency $\omega$ of the perturbation $\Phi(t)$ and the cut-off parameter $N$. It follows from calculations that higher values of $\omega$ are preferable. Numerical experiments show that the value of $C(\mathbf{A}_N,\alpha)$ increases  when $N$ grows and almost does not depend on $\alpha$ for $\omega=8$ and decreases when $\alpha$ grows for $\omega=1$. On the other hand, Table 2 shows that lower values of $\alpha$ correspond to smaller discrepancy $\eta$. This is the reason why the value of $\alpha=0$ has been set in the experiments described below.

Results shown in Table 1 confirm also the Remark made in previous Section. Values of $|H(0)|$ for $\Phi(t)=\sin(8t+\beta_1)\exp(-0.2 t)$ and $\Phi(t)=\sin(t+\beta_2)\exp(-0.2t)$ are $64.02$ and $1.02$ respectively. It is seen from Table 1 that the function $\Phi(t)$ with bigger $|H(0)|=|\Phi''(0)|$ is preferable.

Further we have checked different values of decay coefficient $\nu=0.2\div 10$ of the function $\Phi(t)=\sin(\omega t+\beta) \exp(-\nu t)$. It turned out that bigger values of $\nu$ are preferable because  they decrease $C(\mathbf{A}_N,\alpha)$. For instance, for the value $\nu=10$ and $N$ changing in the range $5\div 20$ the computed values of $C(\mathbf{A}_N,\alpha)$ monotonously raise in the intervals $1.03\div1.6$  and  $1.0\div1.11$ for $\omega=8$ and $\omega=1$ respectively.

In order to obtain admissible values of parameter $N$ for different noise level $\gamma$, we generate synthetic data for $T=12\cdot10^{-9}$ sec, $c=1.5 \cdot 10^8$ m/sec, $l=0.9$ m, $F(x)=\exp(-((x-0.3l)/0.15l)^2)+\exp(-((x-0.7l)/0.1l)^2)$ with function $H(t)=\Phi''(t)$, $\Phi(t)=\sin(8t+1.546)\exp(-0.2 t)$. Different values of $N$ has been tested and the most favorable ones are established. The results are collected in Table 3.

We see in Table 3 that the discrepancy $\eta$ is above the absolute noise level in $g(t)$: $\eta\approx\gamma_1$. This verifies the choice of the values of $N$ and $\alpha=0$ made in the table for the each relative noise level $\gamma$.
\\

{\small \noindent\textbf{Table 3.} Admissible values of the cut-off parameter $N$, corresponding recovery errors $\varepsilon_F$, discrepancy values $\eta$ for  different relative ($\gamma$) and absolute ($\gamma_1=\|\delta g(t)\|_{L_2[0,T]}$) noise levels:
\begin{center}
\begin{tabular}{|c|c|c|c|c|c|c|c|c|}
\hline
$\omega=8$ & $\gamma$ & $ 0\%$ & $ 1\%$ & $ 3\%$ & $5\%$ & $7\%$ & $10\%$ &$20\%$
\\ \cline{2-9} & $\gamma_1$ & $0$ & $0.0064$&$0.019$&$0.032$&$0.045$&$0.064$&$0.128$
\\ \cline{2-9} & $N$ & $20$ & $17$ & $14$ & $11$ & $ 11$ & $ 11 $ &$ 9 $
\\ \cline{2-9} & $\varepsilon_F$ &$0.46\%$ & $0.7\%$ & $1.5\%$ & $2.3\%$ & $ 3\%$ & $4\%$ & $ 7.6\%$
\\ \cline{2-9} & $\eta$ & $0.003$ & $0.008$ & $0.021$ & $0.033$ & $0.044$ & $0.061$ & $0.122$
\\ \hline
$\omega=1$ & $\gamma$ & $ 0\%$ & $ 1\%$ & $ 3\%$ & $5\%$ & $7\%$ & $10\%$ &$20\%$
\\ \cline{2-9} & $\gamma_1$ & $0$ &$0.0076$&$0.023$&$0.038$&$0.053$&$0.076$ &$0.152$
\\ \cline{2-9}& $N$ & $ 20$ & $13$ & $11$ & $10$ & $10$ & $9$ & $9$
\\ \cline{2-9} & $\varepsilon_F$ & $0.6\%$&$2.3\%$ & $3.7\%$ & $4\%$ & $5\%$ & $6.5\%$ & $12\%$
\\ \cline{2-9} & $\eta$ & $0.0005$ & $0.008$ & $0.025$ & $0.036$ & $0.05$ & $0.072$ & $0.143$
\\ \hline
\end{tabular}
\end{center}}

Results of recovery based on parameters and other inputs taken from Table 3 are presented in Fig.1 for the case of $H(t)=\Phi''(t), ~\Phi(t)=\sin(t+1.373)\exp(-0.2t)$.
%
\begin{figure}[!ht]
\begin{center}
\epsfig{file=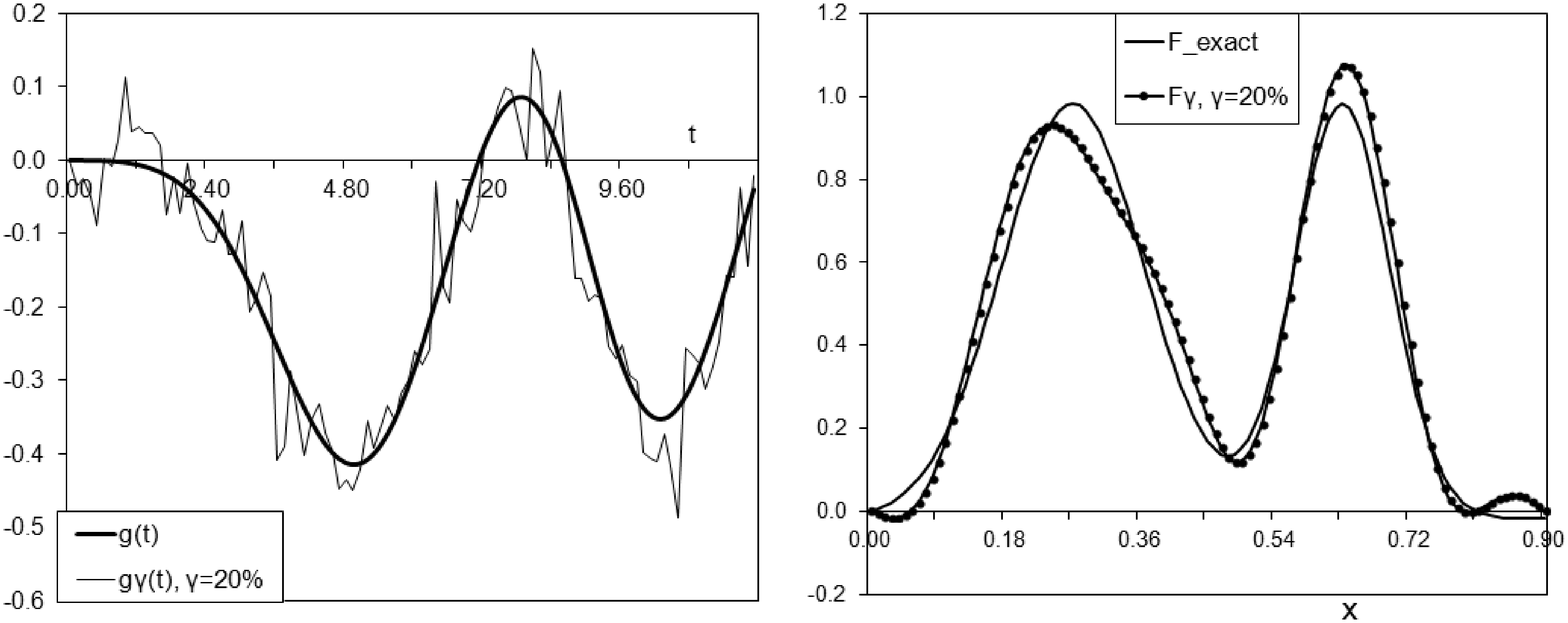,height=6.0cm,width=\textwidth}
\caption{The identified spacewise source $F(x)$ (right figure) from $20\%$ noisy data (left figure)
with parameter $N=9$ defined in Table 3 for $\omega=1$.}
\end{center}
\end{figure}
\begin{figure}[!ht]
\begin{center}
\epsfig{file=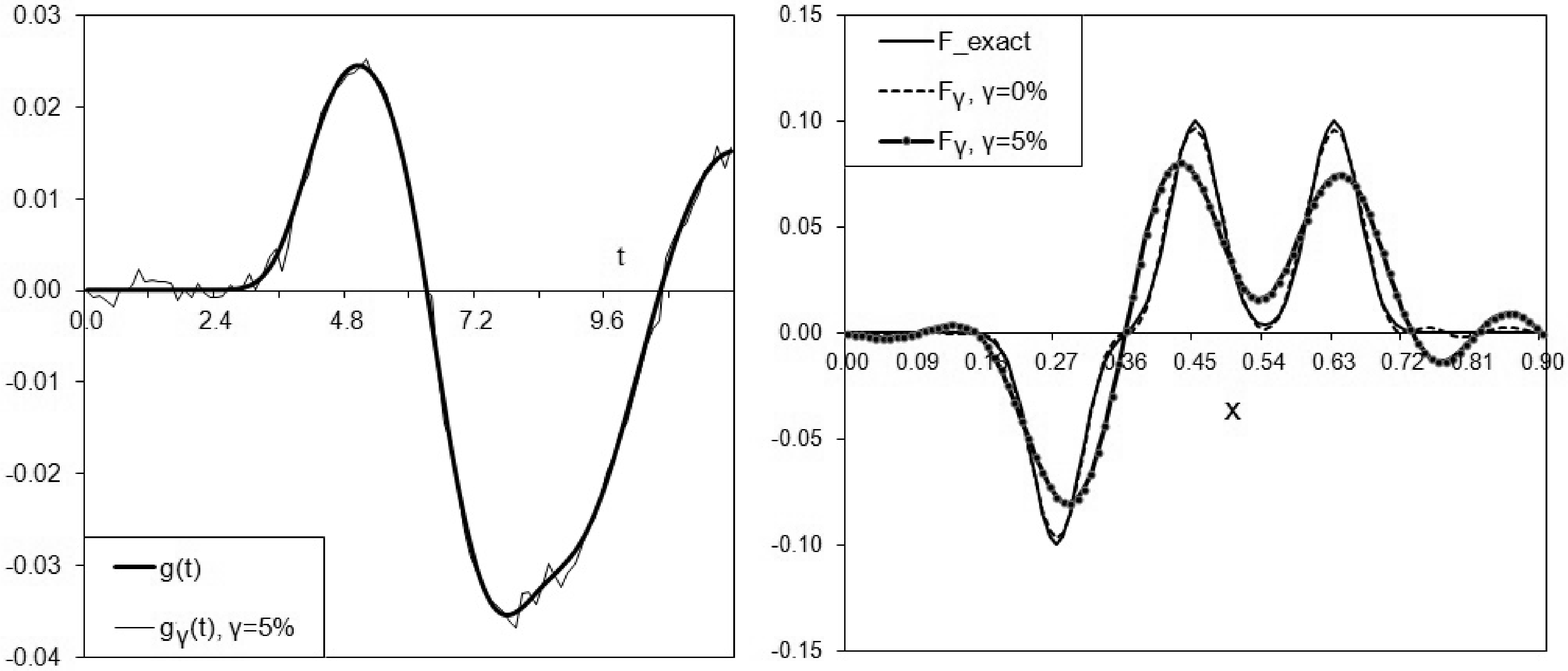,height=6.0cm,width=\textwidth}
\caption{The identified spacewise source $F(x)$ (right figure) from noise free and $5\%$ noisy data (left figure)
with $N=20$ and $10$ for $\Phi(t)=\sin(t+1.373)\exp(-0.2 t)$.}
\end{center}
\end{figure}
\begin{figure}[!ht]
\begin{center}
\epsfig{file=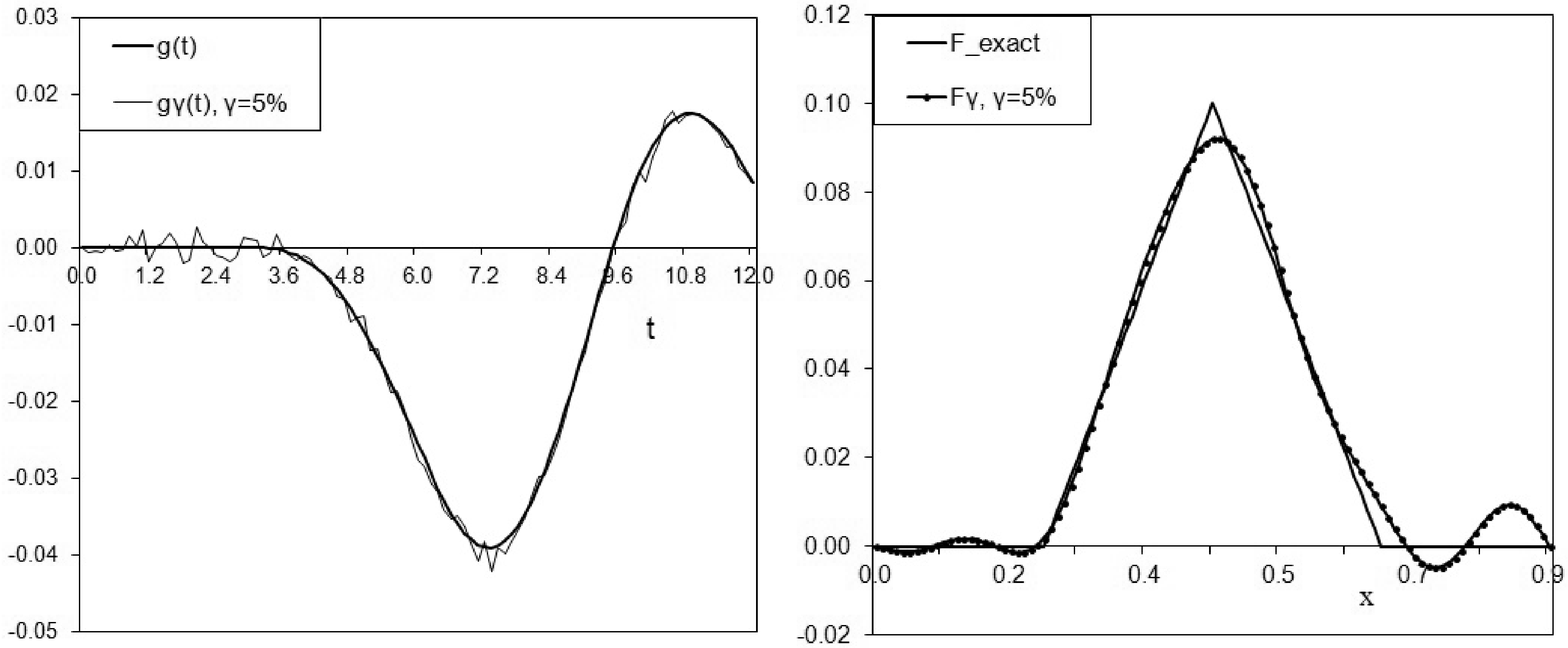,height=6.0cm,width=\textwidth}
\caption{The identified spacewise source $F(x)=\eta(4(x-0.5l)/l)$ (right figure) from $5\%$ noisy data (left figure)
with $N=10$ for $\Phi(t)=\sin(t+1.373)\exp(-0.2 t)$.}
\end{center}
\end{figure}
\begin{figure}[!ht]
\begin{center}
\epsfig{file=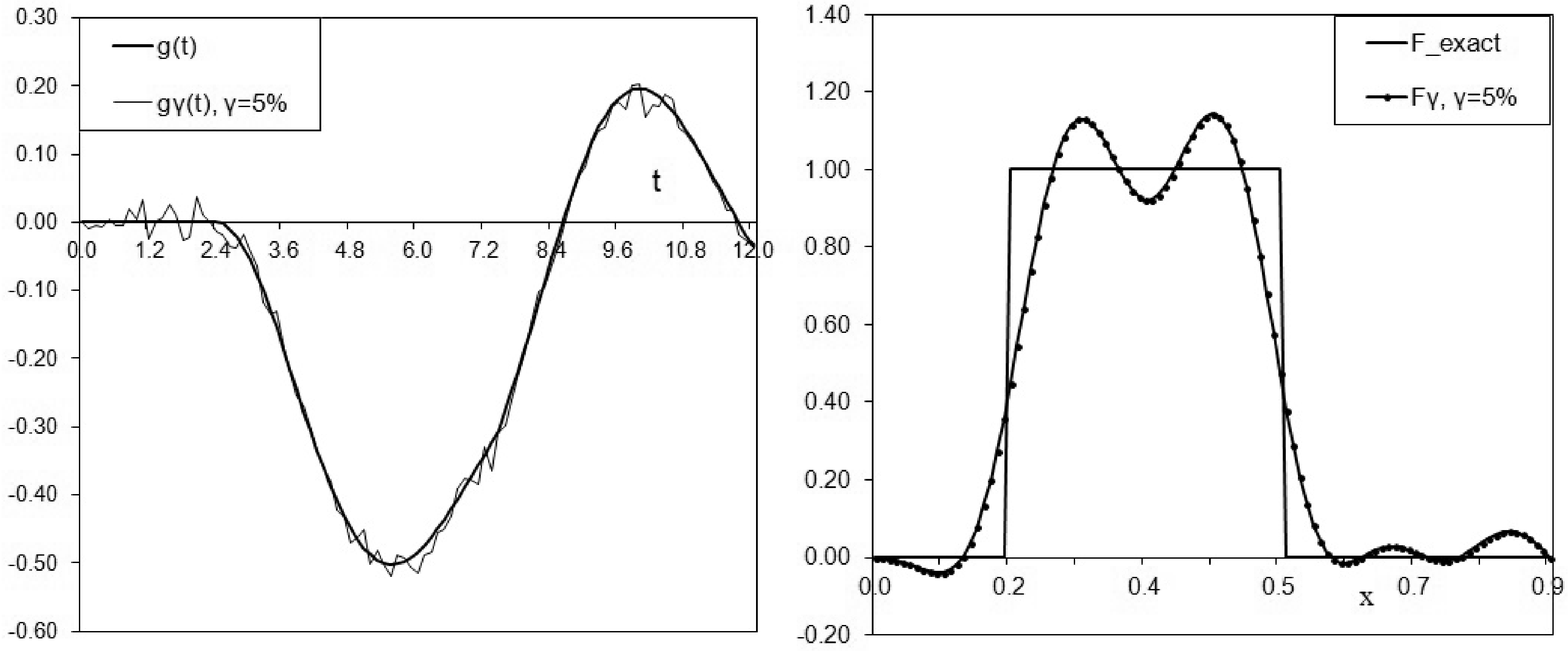,height=6.0cm,width=\textwidth}
\caption{Recovery of  discontinuous source $F(x)$ (right figure) from $5\%$ noisy data (left figure) with parameters listed in Table 2 for $\omega=1$.}
\end{center}
\end{figure}
Since the matrix $\mathbf{A}$ does not depend on $F(x)$, the parameters $N$, $\alpha$ can be defined once for given $H(t)$ and then used again for wide range of functions $F(x)$. Further we take parameters listed in Table 3 to recover other functions, including the combination of three Gaussians $F(x)=-0.1\exp(-((x-0.3l)/0.05l)^2)+0.1\exp(-((x-0.5l)/0.05l)^2)+\exp(-((x-0.7l)/0.05l)^2)$ and the function generated from the standard linear finite element $F(x)=\eta(4(x-0.5l)/l)$. The results are depicted in Fig. 2-3 for the case of $5\%$ noisy data. Then the algorithm has been applied for recovery of discontinuous functions $F(x)$. In that case only approximate agreement  has been achieved. The result is represented in Fig.4.  It follows from numerical simulations that better results are obtained for higher frequency $\omega$ and higher decay coefficient $\nu$  of the function $\Phi(t)=\sin(\omega t+\beta)\exp(-\nu t)$. Numerical simulations show that in the case of smooth $F(x)$ for given parameters $c,T$ and the function $H(t)$, the admissible values of $N$ can be found via numerical experiments on synthetic data.

\section*{\large Conclusion}

In this paper, we studied an inverse problem of identifying the unknown spacewise dependent source $F(x)$ in the one-dimensional wave equation  \mbox{$u_{tt}=c^2u_{xx}+F(x)H(t-x/c)$}, $(x,t)\in \Omega_T$,
which is  treated  as an approximate model of GPR data interpretation process. The case of boundary measured data $g(t):=u(0,t)$ is considered. Perturbation of the media via the radar signal is formulated in terms of the function $H(t-x/c)$ and the non-homogeneity of electrical permittivity is expressed via the function $F(x)$ with finite support in $(0,l)$,  $l=cT/2$. We develop a simple algorithm for reconstruction of a spacewise dependent source term $F(x)$, based on integral formula for the solution of the direct problem with subsequent minimization of the regularized  Tikhonov functional. The proposed algorithm allows one to reconstruct the unknown source from random noisy data up to $10\%$ noise level for a reasonable choice of the function $H(t)$. Note that this method can also be applied to obtain an initial iteration for Conjugate Gradient Algorithm solving the coefficient inverse problem for a hyperbolic equation  $u_{tt}=c^2(x) u_{xx}$ with variable wave propagation speed.

\section*{\large Acknowledgement}
The work of the first author was  supported by the Ministry of Education and Science of Republic of Kazakhstan,  under the Grant No. 316 (13 May, 2016).


\begin{flushleft}
Balgaisha Mukanova,\\
L.N. Gumilyov Eurasian National University,\\
2,Satpayev Str., 010008 Astana, Republic of Kazakhstan,\\
Email: {\tt mukanova\_bg@enu.kz}
\end{flushleft}

\begin{flushleft}
Vladimir G. Romanov,\\
Sobolev Institute of Mathematics,  \\
Novosibirsk 630090, Koptyug prosp., 4, Russia,\\
Email: {\tt romanov@math.nsc.ru}\\

\end{flushleft}

\vspace{1cm}
Received 12.07.2016,   $\> \> \> $ Accepted 05.08.2016

\end{document}